\documentclass{article}
\usepackage{amsmath,amssymb,graphicx,color,url,mathptmx}
\newtheorem{thm}{Theorem}[section]
\newtheorem{dfn}[thm]{Definition}
\newtheorem{prop}[thm]{Proposition}
\newtheorem{cor}[thm]{Corollary}
\newtheorem{lem}[thm]{Lemma}
\newtheorem{rem}[thm]{Remark}

\numberwithin{equation}{section}
\begin{document}

\title{Maximum and minimum of local times for two-dimensional random walk
\footnote{Research partially supported by JSPS Research Fellowships.}}
\author{Yoshihiro Abe \footnote{Kyoto University, Kyoto 606-8502, Japan; 
yosihiro@kurims.kyoto-u.ac.jp}}
\date{}
\maketitle

\begin{abstract}
We obtain the leading orders of the maximum and the minimum of local times for
the simple random walk on the two-dimensional torus at time proportional to the cover time. We also estimate the
number of points with large (or small) values of the local times.
These are analogues of estimates on the two-dimensional 
Gaussian free fields by Bolthausen, Deuschel, and Giacomin 
[Ann. Probab., {\bf 29} (2001)] and Daviaud [Ann. Probab., {\bf 34} (2006)],
but we have different exponents from the case of the Gaussian free field.
\end{abstract}
\textit{MSC 2010}: 60J55; 60J10; 60G70\\
\textit{Keywords}: Local times; Gaussian free fields; Two-dimensional random walks.\\

\section{Introduction} \label{intro}
The theory of local times of random walks is very profound.
It is well-known that local times of random walks have close relationships
with the Gaussian free field(GFF).
The connection goes back to \cite{Dy}.
Eisenbaum, Kaspi, Marcus, Rosen, and Shi \cite{EKMRS}
gave a powerful
equivalence in law called ``generalized second Ray-Knight theorem"
(see Remark \ref{local-time-gff}). 
Using the theorem, Ding, Lee, and Peres \cite{DLP}
established a useful connection between the expected maximum of the GFF and the cover time,
and quite recently
Zhai \cite{Z} strengthened the result
by constructing a coupling of the occupation time filed 
and the GFF (see Theorem \ref{Zhai's thm}).

Much efforts have been made to study local times of the simple random walk
on $\mathbb{Z}^2$.
Erd\H{o}s and Taylor \cite{ET} obtained an estimate on
the maximum of local times of the simple random walk on $\mathbb{Z}^2$
by time $n$.
Dembo, Peres, Rosen and Zeitouni \cite{DPRZ2} improved the result;
they gave the leading order of the maximum 
and estimated the number of ``favorite points"
(see also \cite{R}).
Okada \cite{O} obtained a corresponding estimate on frequently visited sites
in the inner boundary of the random walk range.
Sznitman \cite{S} studied convergences of occupation time fields 
and related the fields to the GFF.

As mentioned above, works \cite{ET, DPRZ2, R} are closely linked to 
the study of extremes of the two-dimensional GFF.
Bolthausen, Deuschel, and Giacomin \cite{BDG}
obtained the leading order of its maximum
(see Remark \ref{remark}).
Daviaud \cite{Da} estimated the number of points with large values of the GFF
(see Remark \ref{remark}).

In this paper, we study the maximum and the minimum of local times
of the simple random walk on the two-dimensional torus at time proportional to the cover time.
While similar work has been done in \cite{BDG, Da} for the GFF,
one cannot apply their results to deduce corresponding local time estimates,
and indeed considerable amounts of efforts are needed to obtain such estimates.
We also note that the exponents for the local times are different from those of the GFF
(see Theorem \ref{maximum} and Remark \ref{remark}).
 
To state our results, we begin with some notation.
We will write $\mathbb{Z}_N^2$ 
to denote the two-dimensional discrete torus with $N^2$ vertices.
Let $X = (X_t)_{t \ge 0}$ be the continuous-time 
simple random walk on $\mathbb{Z}_N^2$
with exponential holding times of parameter $1$.
Let $P_x$ be the law of $X$ starting from $x \in \mathbb{Z}_N^2.$
We define the local time of $X$ by
$$L_t^N (x) := \int_0^t 1_{\{X_s = x \}} ds,~x \in \mathbb{Z}_N^2, t \ge 0,$$
and the inverse local time by 
$$\tau_t := \inf \{s \ge 0 : L_s^N (0) > t \},~t \ge 0.$$
We will take the following time parameter
$$t_{\theta} = t_{\theta}(N) := \frac{4}{\pi} \theta (\log N)^2, \theta > 0.$$
Note that $\tau_{t_{\theta}}$ is approximated by 
$\theta \cdot \frac{4}{\pi} N^2 (\log N)^2$
and that $\frac{4}{\pi} N^2 (\log N)^2$ is close to
the cover time of $\mathbb{Z}_N^2$
(see Lemma \ref{concentration-estimate-inverse-local-time} 
and Theorem 1.1 in \cite{DPRZ3}).
We define sets of ``thick points" and ``thin points" by
\begin{equation} \label{thick-point}
\mathcal{L}_N^{+} (\eta, \theta) := \{x \in \mathbb{Z}_N^2
 : \frac{L_{\tau_{t_{\theta}}}^N (x) - t_{\theta}}
{\sqrt{2t_{\theta}}} \ge \eta \cdot 2\sqrt{2/\pi} \log N \}, ~~
\eta, \theta > 0, 
\end{equation}
\begin{equation} \label{thin-point}
\mathcal{L}_N^{-} (\eta, \theta) := \{x \in \mathbb{Z}_N^2
 : \frac{L_{\tau_{t_{\theta}}}^N (x) - t_{\theta}}
{\sqrt{2t_{\theta}}} \le -\eta \cdot 2\sqrt{2/\pi} \log N \}, ~~
\eta, \theta > 0.
\end{equation}
We will say that $(h_x^N)_{x \in \mathbb{Z}_N^2}$ is the GFF on
$\mathbb{Z}_N^2$ if $(h_x^N)_{x \in \mathbb{Z}_N^2}$ is a centered Gaussian field
with $h_0^N = 0$ and
$\mathbb{E} [h_x^N h_y^N] = E_x [L_{T_0} (y)]$ for all $x, y \in \mathbb{Z}_N^2,$
where $T_0 := \inf \{t \ge 0 : X_t = 0 \}.$
\begin{rem} \label{local-time-gff}
Let $(h_x^N)$ be the GFF on $\mathbb{Z}_N^2$ with a measure $\mathbb{P}$.
The generalized second Ray-Knight theorem \cite{EKMRS}
says that for all $t \ge 0$,
under the measure $P_0 \times \mathbb{P}$,
\begin{equation} \label{second-ray-knight-thm}
\{L_{\tau_t}^N (x) + \frac{1}{2} (h_x^N)^2 : x \in \mathbb{Z}_N^2 \}
= \{\frac{1}{2} (h_x^N + \sqrt{2t})^2 : x \in \mathbb{Z}_N^2 \}~\text{in law}.
\end{equation}
In particular, fixing $N$, we have
\begin{equation} \label{clt}
(\frac{L_{\tau_t}^N (x) - t}{\sqrt{2t}})_{x \in \mathbb{Z}_N^2}
\to (h_x^N)_{x \in \mathbb{Z}_N^2}~\text{in law as}~t \to \infty.
\end{equation}
By (\ref{clt}), one can expect that (\ref{thick-point}) and (\ref{thin-point}) 
will be close
in law to corresponding level sets of the GFF 
(but not exactly, as we see in Theorem \ref{maximum} and Remark \ref{remark}). 
We note that one cannot deduce local time estimates corresponding to
\cite{BDG, Da} from (\ref{second-ray-knight-thm}) or (\ref{clt}).
\end{rem}

We say that a sequence of events $A_N$ holds with high probability if
$\lim_{N \to \infty} P(A_N) = 1.$
We write $|B|$ to denote the cardinality of $B \subset \mathbb{Z}_N^2.$
We now state our results.
\begin{thm} \label{maximum}
(i) For all $\theta> 0$, $\varepsilon > 0,$ 
and $\eta \in (0, 1+\frac{1}{2\sqrt{\theta}})$, the following holds 
with high probability (under $P_0$):
$$N^{2-2(\sqrt{\theta + 2 \eta \sqrt{\theta}} - \sqrt{\theta})^2 - \varepsilon} \le |
\mathcal{L}_N^{+} (\eta, \theta)|
\le N^{2-2(\sqrt{\theta + 2 \eta \sqrt{\theta}} - \sqrt{\theta})^2 + \varepsilon}.$$
Furthermore, for all $\theta > 0$ and $\eta > 1 + \frac{1}{2\sqrt{\theta}}$,
$$|\mathcal{L}_N^{+} (\eta, \theta)| = 0,~~\text{with high probability (under $P_0$)}.$$
(ii) For all $\theta> 1,$ $\varepsilon > 0,$
and $\eta \in (0, 1- \frac{1}{2\sqrt{\theta}})$, the following holds 
with high probability (under $P_0$):
$$N^{2-2(\sqrt{\theta} - \sqrt{\theta - 2 \eta \sqrt{\theta}})^2 - \varepsilon} \le |\mathcal{L}_N^{-} (\eta, \theta)|
\le N^{2-2(\sqrt{\theta} - \sqrt{\theta - 2 \eta \sqrt{\theta}})^2 + \varepsilon}.$$
Furthermore, for all $\theta > 1$ and $\eta > 1 - \frac{1}{2\sqrt{\theta}}$,
$$|\mathcal{L}_N^{-} (\eta, \theta)| = 0,~~\text{with high probability (under $P_0$)}.$$
\end{thm}
The next corollary follows immediately from Theorem \ref{maximum}.
\begin{cor} \label{cor-maximum}
(i) For all $\theta> 0$ and $\varepsilon > 0$, 
the following holds with high probability (under $P_0$):
$$(1 + \frac{1}{2\sqrt{\theta}} - \varepsilon) 2 \sqrt{2/\pi} \log N
\le \frac{\max_{x \in \mathbb{Z}_N^2} 
L_{\tau_{t_{\theta}}}^N (x) - t_{\theta}}
{\sqrt{2t_{\theta}}} 
\le (1 + \frac{1}{2\sqrt{\theta}} + \varepsilon) 2 \sqrt{2/\pi} \log N.$$
(ii) For all $\theta> 1$ and $\varepsilon > 0$, 
the following holds with high probability (under $P_0$):
$$-(1 - \frac{1}{2\sqrt{\theta}} + \varepsilon) 2 \sqrt{2/\pi} \log N
\le \frac{\min_{x \in \mathbb{Z}_N^2} 
L_{\tau_{t_{\theta}}}^N (x) - t_{\theta}}
{\sqrt{2t_{\theta}}} 
\le -(1 - \frac{1}{2\sqrt{\theta}} - \varepsilon) 2 \sqrt{2/\pi} \log N.$$
\end{cor}
\begin{rem} \label{remark}
Set $V_N := [1, N]^2 \cap \mathbb{Z}^2.$
Let $(\tilde{h}_x^N)_{x \in V_N}$
be the GFF on $V_N$ with zero boundary conditions. 
Bolthausen, Deuschel, and Giacomin \cite{BDG}
obtained the leading order of $\max_{x \in V_N} \tilde{h}_x^N$:
for all $\varepsilon > 0$,
$$(1 - \varepsilon) 2 \sqrt{2/\pi} \log N
\le \max_{x \in V_N} \tilde{h}_x^N
\le (1 + \varepsilon) 2 \sqrt{2/\pi} \log N
~~\text{with high probability}.$$
Daviaud \cite{Da} showed that the following holds
with high probability
for all $\varepsilon > 0$ and $\eta \in (0, 1)$:
\begin{equation} \label{daviaud's-result}
N^{2(1 - \eta^2) - \varepsilon} \le 
|\{x \in V_N : \tilde{h}_x^N \ge \eta \cdot 2\sqrt{2/\pi} \log N \} |
\le N^{2(1 - \eta^2) + \varepsilon}.
\end{equation}
We note that one can obtain an estimate similar to (\ref{daviaud's-result}) for the GFF
with periodic boundary conditions by using Theorem \ref{maximum} and Theorem \ref{late-point-thm},
\ref{Zhai's thm} below.
\end{rem}
\begin{rem}
As mentioned before, $t_1$ is close to the cover time for the walk $X$ (see Theorem \ref{late-point-thm} below). 
Thus, it is clear that for $\theta \in (0, 1),$ we have $\min_{x \in \mathbb{Z}_N^2} L_{\tau_{t_{\theta}}}^N (x) = 0.$
\end{rem}

In order to give an intuitive explanation of the exponent in Theorem \ref{maximum}(i),
let us give additional notation.
Let $d(\cdot, \cdot)$ be the $\ell^2$-distance in $\mathbb{Z}_N^2$.
Set $D(x, r) := \{y \in \mathbb{Z}_N^2 : d(x, y) < r \}.$
Fix a subset $A \subset \mathbb{Z}_N^2.$ 
We define its boundary by $\partial A := \{y \in \mathbb{Z}_N^2 : 
y \in \mathbb{Z}_N^2 \backslash A, d(x, y) = 1 ~\text{for some}~x \in A \},$
and the hitting time of $A$ by 
$T_A := \inf \{t \ge 0 : X_t \in A \}.$
We will write $T_x$ to denote $T_{\{x \}}$ for $x \in \mathbb{Z}_N^2.$
Set
$G_A (x, y) := E_x [L_{T_{\partial A}}^N (y)], ~x, y \in \mathbb{Z}_N^2.$
Fix $x \in \mathbb{Z}_N^2$ and $0 < r < R < \frac{N}{2}.$
We define a sequence of stopping times as follows:
$$\tau_x^{(0)} [r, R] := \inf \{t \ge 0 : X_t \in \partial D(x, r) \},$$
$$\sigma_x^{(j)} [r, R] := \inf \{t \ge 0 : X_t \circ 
\theta_{\sum_{i = 0}^{j-1} \tau_x^{(i)} [r, R]} \in \partial D(x, R) \}, ~j \ge 1,$$
$$\tau_x^{(j)} [r, R] := \inf \{t > \sigma_x^{(j)} [r, R]  : X_t \circ 
\theta_{\sum_{i = 0}^{j-1} \tau_x^{(i)} [r, R]} \in \partial D(x, r) \},~j \ge 1,$$
where $\theta_t, t \ge 0$ is the shift operator. 
We define local times of excursions as follows:
$$L_x^{(j)} [r, R] := L_{T_{\partial D(x, R)}}^N (x) \circ \theta_{\sum_{i = 0}^{j-1} 
\tau_x^{(i)} [r, R]},~j \ge 1.$$
We now give heuristics about the exponent in Theorem \ref{maximum}(i).
Let $K_n := n^b e^n n^{3n},$ where $b$ is a positive constant.
We will consider the simple random walk on $\mathbb{Z}_{K_n}^2.$
Set $r_{n, k} := e^n n^{3(n-k)}, k = 0, \dotsc, n$.
For $x \in \mathbb{Z}_{K_n}^2$ and $0 \le \ell \le n-1,$ we write 
$N_{\ell}^x$ to denote the number of excursions from $\partial D(x, r_{n, \ell + 1})$
to $\partial D(x, r_{n, \ell})$ up to time $\tau_{t_{\theta}}$.
By concentration estimates 
(see Lemma \ref{excursion-length} 
and \ref{concentration-estimate-inverse-local-time}),
$$(K_n)^2 t_{\theta} \approx \tau_{t_{\theta}} 
\approx \sum_{j = 0}^{N_0^x} \tau_x^{(j)} [r_{n, 1}, r_{n, 0}]
\approx \frac{2}{\pi} (K_n)^2 \log (\frac{r_{n, 0}}{r_{n, 1}}) N_0^x.$$ 
Thus, we have 
\begin{equation} \label{eq-insight1}
N_0^x \approx 6 \theta n^2 \log n.
\end{equation}
By the law of large numbers, if
$\frac{L_{\tau_{t_{\theta}}}^N (x) - t_{\theta}}{\sqrt{2t_{\theta}}}
\approx \eta \cdot 2 \sqrt{\frac{2}{\pi}} \log K_n$, then we have
\begin{align*}
(\theta + 2 \eta \sqrt{\theta}) \frac{4}{\pi} (\log K_n)^2
&\approx \sum_{j = 1}^{N_{n-1}^x} L_x^{(j)} [r_{n, n}, r_{n, n-1}] \\
&\approx N_{n-1}^x \cdot G_{D(x, r_{n, n-1})} (y, x)
\approx N_{n-1}^x \cdot \frac{2}{\pi} \log (\frac{r_{n, n-1}}{r_{n, n}}),
\end{align*}
where $y$ is a fixed point in $\partial D(x, r_{n, n})$,
and we have used an estimate on Green's functions
(see Lemma \ref{green-function}).
Hence, we have
\begin{equation} \label{eq-insight2}
N_{n-1}^x \approx 6 (\theta + 2 \eta \sqrt{\theta}) n^2 \log n.
\end{equation}
To obtain the order of $|\mathcal{L}_{K_n}^{+} (\eta, \theta)|$,
we should estimate the probability $P_0 (N_{n-1}^x \approx 6 (\theta + 2 \eta \sqrt{\theta}) n^2 \log n)$.
Since for all $1 \le \ell \le n-1$ and $y \in \partial D(x, r_{\ell})$, we have
$P_y (T_{\partial D(x, r_{n, \ell - 1})} < T_{\partial D(x, r_{n, \ell + 1})}) \approx
\frac{1}{2}$ (see Lemma \ref{green-function}),
we can reduce the problem to the case of the simple random walk on $\{0, \dotsc, n \}$;
we need to know the probability of the event that
the walk traverses $6 (\theta + 2 \eta \sqrt{\theta}) n^2 \log n$ times from $n$ to $n-1$
until it crosses $N_0^x$ times from $1$ to $0$.
By this observation, (\ref{eq-insight1}) and a large deviation estimate, we have
$$P_0 (N_{n-1}^x \approx 6 (\theta + 2 \eta \sqrt{\theta}) n^2 \log n)
\approx (K_n)^{- 2(\sqrt{\theta + 2 \eta \sqrt{\theta}} - \sqrt{\theta})^2}.$$
Therefore, if $|\mathcal{L}_{K_n}^{+} (\eta, \theta)|$ is concentrated around
its expectation, we have
$$|\mathcal{L}_{K_n}^{+} (\eta, \theta)| \approx
(K_n)^{2 - 2(\sqrt{\theta + 2 \eta \sqrt{\theta}} - \sqrt{\theta})^2}.$$ \\
The organization of the paper is the following.
Section \ref{preliminary-lem} gives preliminary lemmas.
In Section \ref{sec-maximum}, we prove Theorem \ref{maximum}(i).
The proof is based on the``refined second moment method" 
in \cite{DPRZ, R}.
In Section \ref{sec-minimum}, we prove Theorem \ref{maximum}(ii).

We will write $c_1, c_2, \dotsc$ to denote positive universal constants
whose values are fixed within each argument.
We use $c_1 (\theta), c_2 (\theta) \dotsc$ for positive constants
which depend only on $\theta$.
Given a sequence $(\varepsilon_N)_{N \ge 0}$, we will write $\varepsilon_N = o(1_N)$
if $\lim_{N \to \infty} \varepsilon_N = 0.$

\section{Preliminary lemmas} \label{preliminary-lem}
In this section, we collect lemmas which are useful in the proof of
Theorem \ref{maximum}.
We will use the following basic estimates on the two-dimensional random walk.
See, for example, Theorem 1.6.6, Proposition 1.6.7, and Exercise 1.6.8 in \cite{L}.
\begin{lem} \label{green-function}
(i) There exist $c_1, c_2 >0$ such that the following hold
for all $0 < R < \frac{N}{2}, x \in \mathbb{Z}_N^2,$ and $x_0 \in D(x, R)$:
$$|G_{D(x, R)} (x, x) - \frac{2}{\pi} \log R| \le c_1,$$
$$|G_{D(x, R)} (x_0, x) - \frac{2}{\pi} \log (\frac{R}{d(x_0, x)}) |
\le c_2(\frac{1}{d(x_0, x)} + \frac{1}{R}).$$
(ii) There exist $c_1, c_2 > 0$ such that for all $0 < r < R < \frac{N}{2}$,
$x, x_0 \in \mathbb{Z}_N^2$ with $r < d(x_0, x) < R,$
$$\frac{\log (R/d(x_0, x)) - c_1/r}{\log (R/r)} 
\le P_{x_0} (T_{\partial D(x, r)} < T_{\partial D(x, R)})
\le \frac{\log (R/d(x_0, x)) + c_2/r}{\log (R/r)}.$$
\end{lem}
The following lemma relates time to the number of excursions.
\begin{lem} \label{excursion-length}
There exist $c_1, c_2, c_3$ such that the following holds for all $r, R$ with
$0 < 2r < R < \frac{N}{2}$, $c_1(\frac{1}{r} + \frac{r}{R}) \le \delta \le c_2,$
$x \in \mathbb{Z}_N^2,$ and $M \in \mathbb{N}$:
$$P_0 (\sum_{j=0}^{M} \tau_x^{(j)} [r, R] 
\ge (1+\delta)\frac{2}{\pi} N^2 \log (R/r) M )
\le \exp (- c_3 \delta^2 \frac{\log (R/r)}{\log (N/r)} M).$$
\end{lem} 
{\it Proof.}
The proof is almost the same as that of Lemma 3.2 in \cite{DPRZ} since
Lemma 3.1 of \cite{DPRZ} holds even for the continuous-time simple random walk.
$\Box$ \\
We will use the following moment estimate on local times.
\begin{lem} \label{moment-estimate-local-time}
Fix $x \in \mathbb{Z}_N^2, 0 < R < \frac{N}{2}$, and $x_0 \in D(x, R)$. \\
For all $\beta > 0$,
\begin{equation} \label{moment-local-time-eq}
E_{x_0} [\exp \{-\frac{\beta}{G_{D(x, R)}(x, x)} L_{T_{\partial D(x, R)}} (x) \}]
= 1 - \frac{G_{D(x, R)} (x_0, x)}{G_{D(x, R)} (x, x)} \cdot \frac{\beta}{1+\beta}.
\end{equation}
\end{lem}
{\it Proof.}
By Kac's moment formula
(see, for example, (4) in \cite{FP}), we have for all $k \in \mathbb{N},$
\begin{equation} \label{Kac's formula}
E_{x_0} [(L_{T_{\partial D(x, R)}} (x))^k]
= k! \cdot G_{D(x, R)}(x_0, x) \cdot (G_{D(x, R)} (x, x))^{k-1}.
\end{equation}
The equation (\ref{moment-local-time-eq})
follows immediately from (\ref{Kac's formula}) for all $0 <\beta < 1$.
Regarding both sides of (\ref{moment-local-time-eq}) as analytic functions of
$\beta$, we can show (\ref{moment-local-time-eq}) even for all $\beta \ge 1$
by the uniqueness theorem of analytic functions.
$\Box$ \\
The following is a special version of Lemma 2.1 in \cite{Di}.
\begin{lem} \label{concentration-estimate-inverse-local-time}
There exists $c_1 > 0$ such that for all $t > 0$ and $\lambda \ge 1$,
$$P_0 [|\tau (t) - t N^2| \ge c_1(\sqrt{\lambda t \log N} + \lambda \log N)N^2
] \le 6 e^{-\frac{\lambda}{16}}.$$
\end{lem}
{\it Proof.}
Note that the definition of the inverse local time in \cite{Di} is slightly different from
ours; it corresponds to $\tau_{4t}$ in our notation. 
Since the effective resistances between vertices
in $\mathbb{Z}_N^2$ are of order $\log N$, the statement follows from Lemma 2.1
of \cite{Di}. $\Box$ \\
The following theorem is about the number of ``late points" of $X$.
\begin{thm} \label{late-point-thm}
For all $\varepsilon > 0$ and $\eta \in (0, 1),$
the following holds with high probability (under $P_0$):
$$N^{2 - 2 \eta - \varepsilon} \le
|\{x \in \mathbb{Z}_N^2 : T_x \ge \eta \cdot \frac{4}{\pi} N^2 (\log N)^2 \}|
\le N^{2 - 2 \eta + \varepsilon}.$$
Furthermore, for all $\eta > 1$,
$$|\{x \in \mathbb{Z}_N^2 : T_x \ge \eta \cdot \frac{4}{\pi} N^2 (\log N)^2 \}| = 0~~\text{with high probability (under $P_0$)}.$$
\end{thm}
{\it Proof.} 
Recall that the holding times
of $X$ are independent exponential variables with mean $1$. 
Thus, it is clear that Theorem \ref{late-point-thm} follows immediately
from Proposition 1.1 in \cite{DPRZ}, Theorem 1.1 in \cite{DPRZ3}, and the law of large numbers for the variables.
$\Box$ \\
The following theorem connects ``thick points", ``thin points" and the GFF.
\begin{thm} (Theorem 3.1, \cite{Z})\label{Zhai's thm} 
Let $(h_x^N)_{x \in \mathbb{Z}_N^2}$ be the GFF on $\mathbb{Z}_N^2$.
For all $t > 0$,
$$\{\sqrt{L_{\tau_t}^N (x)} : x \in \mathbb{Z}_N^2 \}
\preceq \frac{1}{\sqrt{2}} \{\max (h_x^N + \sqrt{2t}, 0) : x \in \mathbb{Z}_N^2 \},$$
where $\preceq$ denotes the stochastic domination.
\end{thm}

\section{Proof of Theorem \ref{maximum}(i)} \label{sec-maximum}
Given Theorem \ref{Zhai's thm}, the upper bound of Theorem \ref{maximum}(i) is easy.\\
{\it Proof of the upper bound of Theorem \ref{maximum}(i).}
Fix $\theta > 0$, $\varepsilon > 0,$ and $\eta > 0.$
Let $(h_x^N)_{x \in \mathbb{Z}_N^2}$ be the GFF on $\mathbb{Z}_N^2.$
We have for all $\lambda > 0$,
\begin{align}
&~~~P_0 (|\{x \in \mathbb{Z}_N^2 : \frac{L_{\tau_{t_{\theta}}}^N (x) - t_{\theta}}
{\sqrt{2t_{\theta}}} \ge \eta \cdot 2 \sqrt{2/\pi} \log N \}| 
\ge \lambda) \notag \\ \label{eq0-upper-maximum}
&\le \mathbb{P} (|\{x \in \mathbb{Z}_N^2 : h_x^N \ge 
(\sqrt{\theta + 2 \eta \sqrt{\theta}} - \sqrt{\theta}) \cdot 2 \sqrt{2/\pi}
\log N \}| 
\ge \lambda) 
\\
\label{eq1-upper-maximum}
&= \mathbb{P} (|\{x \in \mathbb{Z}_N^2 : h_x^N 
\le - (\sqrt{\theta + 2 \eta \sqrt{\theta}} - \sqrt{\theta}) \cdot 2 \sqrt{2/\pi} 
\log N \}| \ge \lambda) \\ \label{eq2-upper-maximum}
&\le P_0 (|\{x \in \mathbb{Z}_N^2 : 
L_{\tau_{t_{(\sqrt{\theta + 2 \eta \sqrt{\theta}} - \sqrt{\theta})^2}}}^N (x) = 0 \}| 
\ge \lambda)  \\ \label{eq3-upper-maximum}
&= P_0 (|\{x \in \mathbb{Z}_N^2 : T_x > 
\tau_{t_{(\sqrt{\theta + 2 \eta \sqrt{\theta}} - \sqrt{\theta})^2}} \}| \ge \lambda), 
\end{align} 
where we have used the symmetry of the GFF in (\ref{eq1-upper-maximum})
and Theorem \ref{Zhai's thm} in (\ref{eq0-upper-maximum}) and
(\ref{eq2-upper-maximum}).
Take $\lambda = N^{2 - 2(\sqrt{\theta + 2 \eta \sqrt{\theta}} - \sqrt{\theta})^2 +
 \varepsilon},$
 $\varepsilon > 0.$
By Lemma \ref{concentration-estimate-inverse-local-time}
and Theorem \ref{late-point-thm},
the probability in (\ref{eq3-upper-maximum}) goes to $0$ as $N \to \infty$.~$\Box$ \\

From now on, we prove the lower bound of Theorem \ref{maximum}(i)
by applying the methods in \cite{DPRZ, R}.
First, we define the notion that a point is ``successful".
Set
\begin{equation} \label{space-parameter}
r_{n, k} := e^n n^{3(n-k)},~~0 \le k \le n, ~K_n := n^{\Bar{\gamma}} r_{n, 0},
\end{equation}
where $\Bar{\gamma} \in [b, b + 4]$ and $b$ is a sufficiently large positive constant.
Since $K_n$'s take values over all sufficiently large positive integers, we may only
 consider the subsequence.
From now, we will consider the simple random walk on $\mathbb{Z}_{K_n}^2$.
Given $\eta \in (0, 1 + \frac{1}{2 \sqrt{\theta}})$, we set
$$n_{\ell} := \lceil 
6 (1-1/n^{1/4}) 
\{ \sqrt{\theta} + (\sqrt{\theta+ 2\eta \sqrt{\theta}} - \sqrt{\theta})  
\frac{\ell}{n}  \}^2
n^2 \log n \rceil,~~0 \le \ell \le n-1.$$
For $x \in \mathbb{Z}_{K_n}^2$ and $1 \le \ell \le n-1$,
\begin{align*}
N_{\ell}^x &~\text{is the number of excursions from}~ \partial D(x, r_{n, \ell + 1})
~\text{to}~ \partial D(x, r_{n, \ell}) \\
&~~~~~\text{up to random time}~ 
\sum_{j = 0}^{n_0} \tau_x^{(j)}[r_{n, 1}, r_{n, 0}]. 
\end{align*}

\begin{figure}
\begin{center}
\includegraphics[width=40.0mm, clip]{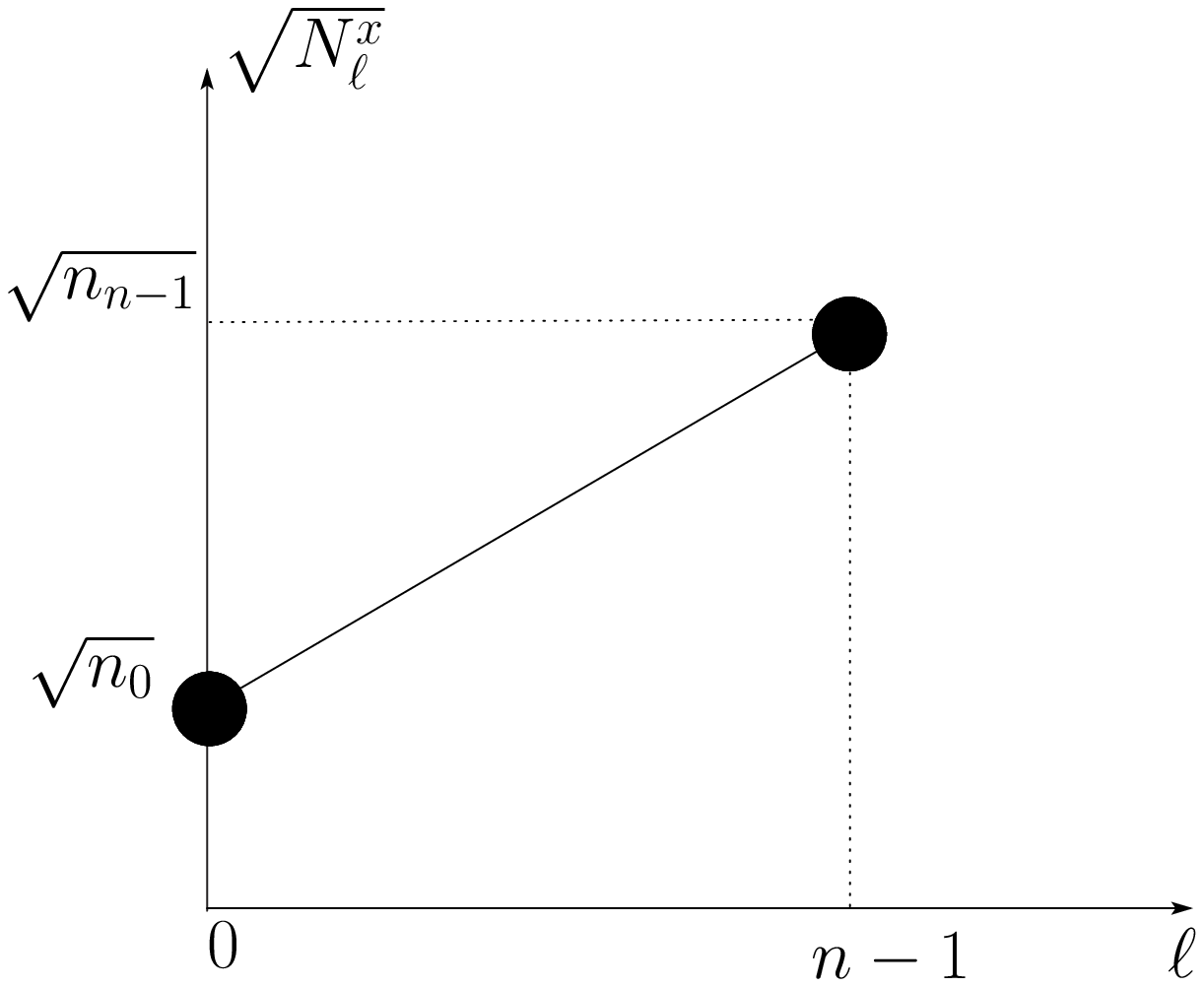}
\end{center}
\caption{If $x$ is successful, $(\sqrt{N_{\ell}^x})_{0 \le \ell \le n-1}$
behaves like a linear function.} \label{fig}
\end{figure}

\begin{dfn} \label{successful}
Fix $x \in \mathbb{Z}_{K_n}^2$. We will say that
$x$ is successful if 
$$|N_{\ell}^x - n_{\ell}| \le n,~~\text{for}~1 \le \ell \le n-1.$$ 
\end{dfn}
\begin{rem}
We give an intuition about Definition \ref{successful}.
Assume that $\frac{L_{\tau_{t_{\theta}}}^N (x) - t_{\theta}}{\sqrt{2t_{\theta}}}
\approx \eta \cdot 2 \sqrt{2/\pi} \log K_n$.
We already know that under this assumption, 
$N_0^x \approx 6 \theta n^2 \log n$ and
$N_{n-1}^x \approx 6 (\theta + 2 \eta \sqrt{\theta}) n^2 \log n$
(recall (\ref{eq-insight1}) and (\ref{eq-insight2})).
Due to a recent work by Belius and Kistler \cite{BK}, 
one expects that
conditioned on
$\sqrt{N_0^x} \approx \sqrt{n_0}$ and $\sqrt{N_{n-1}^x} \approx \sqrt{n_{n-1}}$, 
$(\sqrt{N_{\ell}^x})_{0 \le \ell \le n-1}$ behaves 
roughly like
a Brownian bridge from $\sqrt{n_0}$ to $\sqrt{n_{n-1}}$.
Therefore, we see that 
$(\sqrt{N_{\ell}^x})_{0 \le \ell \le n-1}$ 
would typically look like a linear function in $\ell$
with $\sqrt{N_0^x} \approx \sqrt{n_0}$ and $\sqrt{N_{n-1}^x} \approx \sqrt{n_{n-1}}$
(see Figure \ref{fig}).
We used this insight in Definition \ref{successful}.
Note that our framework is quite different from those in \cite{DPRZ, R} and
so is the definition of ``successful".
\end{rem}
The lower bound of Theorem \ref{maximum}(i)
follows from the following three propositions.
\begin{prop} \label{localtime_excursion}
For all $\theta > 0$ and $\eta \in (0, 1+\frac{1}{2\sqrt{\theta}})$,
the following holds with high probability (under $P_0$):
$$\{x \in \mathbb{Z}_{K_n}^2 \backslash D(0, r_{n, 0}) : x~\text{is successful} \}
\subset \mathcal{L}_{K_n}^{+}
( (1-1/\log \log n)(1-1/n)^2\eta, \theta).$$
\end{prop} 
\begin{prop} \label{probab_successful}
For all $\theta > 0$, $\eta \in (0, 1+\frac{1}{2\sqrt{\theta}})$, and
$x \in \mathbb{Z}_{K_n}^2 \backslash D(0, r_{n, 0})$,
$$P_0 (x~\text{is successful}) = (1 + o(1_n))q_n,$$
where $q_n$ satisfies the following: there exists 
$c_1(\theta), c_2 (\theta)> 0$ such that
$$e^{- c_1 (\theta) n \log \log n} (K_n)^{-2(\sqrt{\theta + 2 \eta \sqrt{\theta}} - 
\sqrt{\theta})^2} \le q_n 
\le e^{- c_2 (\theta) n \log \log n} (K_n)^{-2(\sqrt{\theta + 2 \eta \sqrt{\theta}} - 
\sqrt{\theta})^2}.$$
\end{prop}
\begin{prop} \label{correlation_successful}
Let $q_n$ be given in Proposition \ref{probab_successful}.
Fix $\theta > 0$ and $\eta \in (0, 1+\frac{1}{2\sqrt{\theta}})$.
For $x, y \in \mathbb{Z}_{K_n}^2$, set
$$\ell(x, y) := \min \{\ell : D(x, r_{n, \ell} +1) \cap D(y, r_{n, \ell} + 1) = \emptyset \} \wedge n.$$
(i) There exist $c_1(\theta), c_2 (\theta)> 0$
such that for all $x, y \in \mathbb{Z}_{K_n}^2
\backslash D(0, r_{n, 0})$ with $1 \le \ell(x, y) \le n-2,$
\begin{align*}
&~~P_0 (x~\text{and}~y~\text{are successful}) \\
&\le n^{c_1 (\theta)} e^{c_2 (\theta) \ell \log \log n} (q_n)^2
\cdot \exp \{6 (\sqrt{\theta + 2\eta \sqrt{\theta}} - \sqrt{\theta})^2 
\ell \log n \}.
\end{align*}
(ii) For all $x, y \in \mathbb{Z}_{K_n}^2 \backslash D(0, r_{n, 0})$ 
with $\ell(x, y) = 0,$
$$P_0 (x~\text{and}~y~\text{are successful}) = (1 + o(1_n)) (q_n)^2.$$
\end{prop}

\textit{Proof of the lower bound of Theorem \ref{maximum}(i) 
via Proposition \ref{localtime_excursion}-\ref{correlation_successful}.}
Fix $\theta > 0$ and $\eta 
\in (0, 1+\frac{1}{2\sqrt{\theta}})$. 
 Set
$$Z_n := \sum_{x \in \mathbb{Z}_{K_n}^2 \backslash D(0, r_{n, 0})} 
1_{\{x~\text{is successful} \}},~
W_{n, \ell} := \sum_{\begin{subarray}{c} 
x, y \in \mathbb{Z}_{K_n}^2 \backslash D(0, r_{n, 0})
\\ \ell(x, y) = \ell \end{subarray}} 1_{\{x~\text{and $y$ are successful} \}}.$$
We have
$E_0 [(Z_n)^2] \le \sum_{\ell = 0}^{n} E_0 [W_{n, \ell}].$
Recall (\ref{space-parameter}). Taking $b$ large enough,
by Proposition \ref{probab_successful} and \ref{correlation_successful}, we have
$$E_0[W_{n, 0}] \le (1 + o(1_n)) (K_n)^4 (q_n)^2,~
\sum_{\ell = 1}^{n} E_0[W_{n, \ell}] \le o(1_n) (K_n)^4 (q_n)^2.$$
Thus, we have
\begin{equation} \label{secondmoment}
E_0[(Z_n)^2] \le (1+o(1_n)) (K_n)^4 (q_n)^2.
\end{equation}
By (\ref{secondmoment}), Proposition \ref{probab_successful}, 
and the Paley-Zygmund inequality, the following holds with high probability:
\begin{equation} \label{number_successful}
Z_n \ge e^{-n(\log \log n)^2} 
(K_n)^{2-2(\sqrt{\theta + 2 \eta \sqrt{\theta}} - \sqrt{\theta})^2 }.
\end{equation}
The lower bound of Theorem \ref{maximum}(i)
follows from (\ref{number_successful})
and Proposition \ref{localtime_excursion}. $\Box$ \\

For the rest of this section, we will prove Proposition \ref{localtime_excursion}-\ref{correlation_successful}.\\
\textit{Proof of Proposition \ref{localtime_excursion}.}
We will prove the following:
\begin{align} \label{aim-localtime-excursion}
&P_0 [\frac{L_{\tau_{t_{\theta}}}^N (x) - t_{\theta}}
{\sqrt{2t_{\theta}}} < (1-1/\log \log n)(1-1/n)^2 
\eta \cdot 2\sqrt{2/\pi} \log K_n, \notag\\
&~~~~~~~~~~~~~~~~~~~~~~~~~~~~~~~\text{and}~x~
\text{is successful for some}~ 
x \in \mathbb{Z}_{K_n}^2 \backslash D(0, r_{n, 0})] \notag\\
&\to 0~\text{as}~n \to \infty. 
\end{align}
The statement in Proposition \ref{localtime_excursion} follows immediately from this.
The probability in (\ref{aim-localtime-excursion}) is bounded above by 
$I_1+ I_2 + I_3,$
where
\begin{align} \label{main-term-localtime-excursion}
&I_1 := P_0 [\text{For some}~ 
x \in \mathbb{Z}_{K_n}^2 \backslash D(0, r_{n, 0}),
\sum_{j=1}^{n_{n-1}-n} L_x^{(j)}[r_{n,n}, r_{n, n-1}] \notag\\
&~~~~~~~~~~~~~~~~~~~~~~~~~~~~~~~~
< (\theta + 2(1-1/\log \log n)(1-1/n)^2
\eta\sqrt{\theta})\frac{4}{\pi}(\log K_n)^2 ],
\end{align}
\begin{equation} \label{second-term-localtime-excursion}
I_2 := P_0 [\sum_{j=0}^{n_0} \tau_x^{(j)}[r_{n,1}, r_{n,0}] > \lambda_n 
~\text{for some}~ x \in \mathbb{Z}_{K_n}^2 \backslash D(0, r_{n, 0})],
\end{equation}
\begin{equation} \label{third-term-localtime-excursion}
I_3 := P_0 (\tau_{t_{\theta}} \le \lambda_n ),
\end{equation}
where $\lambda_n := (1+1/n^{1/4}) 2/\pi (K_n)^2 
\log (r_{n,0}/r_{n,1}) n_0 ~~(\le (1-1/\sqrt{n}) (K_n)^2 t_{\theta}).$
By Lemma \ref{excursion-length} and 
\ref{concentration-estimate-inverse-local-time}, we have
$I_2 = o(1_n), I_3 = o(1_n).$ \\
From now, we will prove $I_1 = o(1_n).$
Fix $x \in \mathbb{Z}^2_{K_n} \backslash D(0, r_{n,0}).$ Set
$$P_x := P_0 [\sum_{j=1}^{n_{n-1}-n} L_x^{(j)}[r_{n,n}, r_{n, n-1}]
< (\theta + 2(1-1/\log \log n)(1-1/n)^2
\eta\sqrt{\theta}) 4/\pi (\log K_n)^2].$$
By Lemma \ref{green-function}(i) and
\ref{moment-estimate-local-time} together with 
the Chebyshev inequality and the strong Markov property, 
we have for $\varphi > 0$,
\begin{align*}
P_x
& \le \exp \{\frac{\varphi}{G_{D(x, r_{n,n-1})} (x,x)} 
(\theta + 2(1-1/\log \log n)(1-1/n)^2 
\eta \sqrt{\theta}) 4/\pi (\log K_n)^2 \} \\
&~~\cdot \{ \max_{y \in \partial D(x, r_{n,n})} E_y [ \exp 
\{- \frac{\varphi}{G_{D(x, r_{n,n-1})} (x, x)} L_{T_{\partial D(x, r_{n,n-1})}} (x)
 \} ] \}^{n_{n-1}-n} \\
&\le e^{c_1 (\theta) (\varphi n\log n + \frac{\varphi}{1+\varphi} (\log n)^3)} 
\exp \{- 18n (\log n)^2 f_n (\varphi) \},
\end{align*}
where
\begin{align*}
f_n (\varphi) &:= - \varphi (\theta + 2(1-1/\log \log n)(1-1/n)^2 \eta\sqrt{\theta}) \\
&+ \frac{\varphi}{1+\varphi}(1- 1/n^{1/4})(1-1/n)^2 (\theta + 2\eta \sqrt{\theta}).
\end{align*}
Taking $\varphi$ at which $f_n(\varphi)$ attains the maximum,
we have
$$P_x \le e^{c_2 (\theta) \frac{n\log n}{\log \log n}} \exp \{- c_3 (\theta) n
(\log n/\log \log n)^2 \} = o ( (K_n)^{-2}).$$
Therefore, we have proved $I_1 = o(1_n)$ and 
(\ref{aim-localtime-excursion}). $\Box$ \\

\textit{Proof of Proposition \ref{probab_successful}.}
Fix $x \in \mathbb{Z}_{K_n}^2 \backslash D(0, r_{n,0}).$
By Lemma \ref{green-function}(ii) and the strong Markov property, we have
\begin{equation} \label{eq-successful}
P_0 (x~\text{is successful}) = 
\sum P(N_{\ell}^x = m_{\ell}~\text{for}~1 \le \ell \le n-1)
= (1 + o(1_n)) q_n,
\end{equation}
where
\begin{equation} \label{successfulprobability}
q_n := \sum \prod_{\ell = 1}^{n-1} \begin{pmatrix} m_{\ell} + m_{\ell-1} -1 \\ m_{\ell}
\end{pmatrix} \Bigl(\frac{1}{2} \Bigr)^{m_{\ell} + m_{\ell-1}}.
\end{equation}
Here the summations in 
(\ref{eq-successful}) and (\ref{successfulprobability}) 
are over all $m_1, \dotsc, m_{n-1}$
with $|m_i - n_i| \le n$ for $1 \le i \le n-1.$
By the Stirling formula, we have for all $m_i$ with $|m_i - n_i| \le n, 1 \le i \le n-1,$
\begin{align} 
&\begin{pmatrix} m_{\ell} + m_{\ell-1} -1 \\ m_{\ell}
\end{pmatrix} \Bigl(\frac{1}{2} \Bigr)^{m_{\ell} + m_{\ell-1}}
\ge \frac{c_1 (\theta)}{\sqrt{m_{\ell}}}
\cdot \frac{(m_{\ell} + m_{\ell-1})^{m_{\ell} + m_{\ell-1}}}
{(m_{\ell})^{m_{\ell}} (m_{\ell-1})^{m_{\ell-1}}} 
\Bigl(\frac{1}{2} \Bigr)^{m_{\ell} + m_{\ell-1}} \notag\\
&\ge c_1 (\theta) \cdot (m_{\ell})^{-1/2} \cdot 
\exp \{m_{\ell-1} f(\frac{m_{\ell}}{m_{\ell-1}}) \} \notag\\
&\ge c_1 (\theta) \cdot (m_{\ell})^{-1/2} \cdot 
\exp \{- m_{\ell-1} (\frac{1}{4} (\frac{m_{\ell}}{m_{\ell-1}} - 1)^2
+ c_2 (\theta)
 |\frac{m_{\ell}}{m_{\ell-1}} - 1|^3)  \} \label{stirling_estimate1} \\
&\ge c_3 (\theta) n^{-1} (\log n)^{-1/2}
\exp \{- 6 (\sqrt{\theta + 2\eta \sqrt{\theta}} - \sqrt{\theta})^2 
\log n \}, \label{stirling_estimate2}
\end{align}
where
$f(u) := (1 + u) \log (1 + u) - u \log u - (1+ u) \log 2, u > 0$
and we have used the Taylor expansion of $f$ around $1$ in (\ref{stirling_estimate1}). 
Therefore, we have by (\ref{successfulprobability}) and (\ref{stirling_estimate2})
\begin{align*}
q_n &\ge n^{n-1} (c_3 (\theta) n^{-1}(\log n)^{-1/2})^{n-1}
\exp \{- 6 (\sqrt{\theta + 2\eta \sqrt{\theta}} - \sqrt{\theta} )^2 
n\log n \} \\
&\ge (c_4 (\theta) (\log n)^{-1/2})^{n-1}
(K_n)^{-2(\sqrt{\theta + 2 \eta \sqrt{\theta}} - \sqrt{\theta})^2}.
\end{align*}
By a similar argument, we can obtain the upper bound of $q_n$. $\Box$ \\

In order to prove Proposition \ref{correlation_successful}, we make some preparations.
Fix $x \in \mathbb{Z}_N^2$ and $0 \le \ell \le n-1.$
Set $e^{(0)} := (X_t : 0 \le t \le T_{\partial D(x, r_{n, \ell + 1})}),$ and
$$e^{(i)} := (X_{t + \sum_{j=0}^{i-1} \tau_x^{(j)} [r_{n, \ell + 1}, r_{n, \ell}]} :
\sigma_x^{(i)} [r_{n, \ell + 1}, r_{n, \ell}] \le t \le 
\tau_x^{(i)}[r_{n, \ell + 1}, r_{n, \ell}] ),~i \ge 1.$$
Let $\mathcal{G}_{\ell}^x := \sigma (e^{(i)} : i \ge 0).$ 
We will use the following lemma iteratively.
\begin{lem} \label{decoupling-lem}
Fix $\eta \in (0, 1 + \frac{1}{2\sqrt{\theta}})$ and $\theta > 0.$
There exists $\varepsilon_n$ with $\lim_{n \to \infty} \varepsilon_n = 0$
such that the following holds for all $0 \le \ell \le n - 2$, $m_{\ell}$ with
$|m_{\ell} - n_{\ell}| \le n, m_{\ell +1}, \dotsc, m_{n-1} > 0$, 
and $x \in \mathbb{Z}_{K_n}^2 \backslash D(0, r_{n, 0}):$
\begin{align*}
&P_0 (N_i^x = m_i~\text{for all}~i = \ell, \dotsc, n-1 | \mathcal{G}_{\ell}^x) \\
&= (1 + \varepsilon_n) 
P_0 (N_i^x = m_i~\text{for all}~i = \ell + 1, \dotsc, n-1 | N_{\ell}^x = m_{\ell}) \cdot
1_{\{N_{\ell}^x = m_{\ell} \}}.
\end{align*}
\end{lem}
{\it Proof.} 
The proof is almost the same as that of Corollary 5.1 in \cite{DPRZ}
since Lemma 2.4 in \cite{DPRZ} holds even for the continuous-time simple
random walk. $\Box$ \\

\textit{Proof of Proposition \ref{correlation_successful}.}
Proposition \ref{correlation_successful}(ii) 
follows immediately from Lemma \ref{decoupling-lem}.
We will prove Proposition \ref{correlation_successful}(i).
Fix $x, y \in \mathbb{Z}_{K_n}^2 \backslash D(0, r_{n,0})$ with $1 \le \ell(x,y) \le n-2.$ 
We will write
$$\ell := \ell(x,y).$$ 
By Lemma \ref{decoupling-lem},
\begin{align}
&~~~P_0 (x~\text{and}~y~\text{are successful}) \notag \\
&\le P_0 (|N_i^x - n_i| \le n~\text{for}~i=1, \dotsc, \ell-3, \ell, \dotsc, n-1, 
\notag \\
&~~~~~~~~~~~~~~~~~~~~~~~~~\text{and}~|N_i^y - n_i| \le n~\text{for}~i= \ell, 
\dotsc, n-1) 
\notag \\
&\le (1+o(1_n)) 
P_0 (|N_i^x - n_i| \le n~\text{for}~i=1, \dotsc, \ell-3, \ell, \dotsc, n-1 ) 
\notag \\
&\cdot \sum_{m_{\ell}:|m_{\ell} - n_{\ell}| \le n}
P_0 (|N_i^y - n_i| \le n~\text{for}~i= \ell + 1, \dotsc, n-1|N_{\ell}^y = m_{\ell}).
 \label{eq-decoupling-1}
\end{align}
We will prove the following:
\begin{align} \label{key1_correlation}
\sum_{m_{\ell}:|m_{\ell} - n_{\ell}| \le n}
&~~P_0 (|N_i^y - n_i| \le n~\text{for}~i= \ell + 1, \dotsc, n-1|N_{\ell}^y = m_{\ell})
\notag \\
&\le c_1 (\theta) n e^{c_2 (\theta) \ell \log \log n} q_n
\exp \{6(\sqrt{\theta + 2\eta \sqrt{\theta}} -\sqrt{\theta})^2 
\ell \log n \},
\end{align}
\begin{align} \label{key2_correlation}
&P_0 (|N_i^x - n_i| \le n~\text{for}~i=1, \dotsc, \ell-3, \ell, \dotsc, n-1) 
\notag \\
&\le c_3 (\theta) n^{c_4 (\theta)} e^{c_5 (\theta) \ell \log \log n} q_n.
\end{align}
Proposition \ref{correlation_successful}(i) follows from (\ref{eq-decoupling-1}),
(\ref{key1_correlation}) and (\ref{key2_correlation}). \\
First, we prove (\ref{key1_correlation}). 
By Proposition \ref{probab_successful} and Lemma \ref{decoupling-lem},
\begin{align} \label{successful_decomposition1}
(1+o(1_n)) &q_n = P_0 (y~\text{is successful}) \notag\\
&\ge (1+o(1_n)) \sum_{m_{\ell}:|m_{\ell} - n_{\ell}| \le n} 
P_0 (|N_i^y - n_i| \le n~\text{for}~
1 \le i \le \ell-1,~~N_{\ell}^y = m_{\ell}) \notag \\
&~~~~~~~~~~~~~~~~~~~~~~~~~~~~
\cdot P_0 (|N_i^y - n_i| \le n~\text{for}~\ell+1 \le i \le n-1 |
N_{\ell}^y = m_{\ell}).
\end{align}
By a similar argument to the proof of Proposition \ref{probab_successful}, we have for all 
$m_{\ell}$ with $|m_{\ell} - n_{\ell}| \le n,$
\begin{align} \label{successful_decomposition2}
&P_0 (|N_i^y - n_i| \le n ~\text{for}~1 \le i \le \ell-1,~N_{\ell}^y = m_{\ell}) 
\notag \\
&\ge n^{-1} e^{-c_6 (\theta) \ell \log \log n}
\exp \{-6 (\sqrt{\theta+2\eta \sqrt{\theta}} -\sqrt{\theta})^2 
\ell \log n \}.
\end{align} 
(\ref{key1_correlation}) follows from (\ref{successful_decomposition1}) and
(\ref{successful_decomposition2}). \\
Next, we prove (\ref{key2_correlation}).
By Lemma \ref{decoupling-lem}, (\ref{key1_correlation}), and a similar argument to
the proof of Proposition \ref{probab_successful}, we have
\begin{align*}
&P_0 (|N_i^x - n_i| \le n~\text{for}~i = 1, \dotsc, \ell-3, \ell, \dotsc, n-1)
\\
&\le (1+o(1_n)) P_0 (|N_i^x - n_i| \le n~\text{for}~
1 \le i \le \ell-3) \\
&\cdot \sum_{m_{\ell}:|m_{\ell} - n_{\ell}| \le n}
P_0 (|N_i^x - n_i| \le n~\text{for}~\ell+1 \le i \le n-1 |
N_{\ell}^x = m_{\ell}) \\
&\le e^{-c_7 (\theta) (\ell-3) \log \log n} 
\exp \{-6 (\sqrt{\theta+2\eta \sqrt{\theta}} -\sqrt{\theta})^2 
(\ell-3) \log n \} \\
&\cdot c_1(\theta) q_n n e^{c_2(\theta) \ell \log \log n} 
\exp \{6 (\sqrt{\theta +2\eta\sqrt{\theta}} - \sqrt{\theta})^2 
\ell \log n \} \\
&\le q_n n^{c_8 (\theta)} e^{c_9 (\theta) \ell \log \log n}.
\end{align*}
Therefore, we have proved (\ref{key2_correlation}). $\Box$

\section{Proof of Theorem \ref{maximum}(ii)} \label{sec-minimum}
In this section, we prove Theorem \ref{maximum}(ii).
First, we show the lower bound. \\
{\it Proof of the lower bound of Theorem \ref{maximum}(ii).}
Fix $\theta > 1$, $\varepsilon >0$, and $\eta \in (0, 1 - \frac{1}{2\sqrt{\theta}}).$
We have for all $\lambda > 0$,
\begin{align}
&~~~P_0 ( |\{x \in \mathbb{Z}_N^2 : 
\frac{L_{\tau_{t_{\theta}}}^N (x) - t_{\theta}}
{\sqrt{2t_{\theta}}} \le -\eta \cdot 2\sqrt{2/\pi} \log N \} |
\ge \lambda) \notag \\ \label{eq1-lower-minimum}
& \ge \mathbb{P} ( |\{x \in \mathbb{Z}_N^2 : h_x^N \le 
-(\sqrt{\theta} - \sqrt{\theta - 2 \eta \sqrt{\theta}}) 2 \sqrt{2/\pi} \log N \}|
\ge \lambda )  \\ \label{eq2-lower-minimum}
& = \mathbb{P} ( |\{x \in \mathbb{Z}_N^2 : h_x^N \ge 
(\sqrt{\theta} - \sqrt{\theta - 2 \eta \sqrt{\theta}}) 2 \sqrt{2/\pi} \log N \}|
\ge \lambda)  \\ 
&= \mathbb{P} (|\{x \in \mathbb{Z}_N^2 : \frac{1}{\sqrt{2}}(h_x^N + \sqrt{2 t_{(\theta - 2 \eta \sqrt{\theta})}}) 
\ge \sqrt{\theta} 
\cdot 2/\sqrt{\pi} 
\log N \}| \ge \lambda)  \notag \\ \label{eq3-lower-minimum}
&\ge P_0 (|\mathcal{L}_N^{+} (\eta\sqrt{\theta}/\sqrt{\theta - 2 \eta \sqrt{\theta}}, \theta - 2 \eta \sqrt{\theta})|
\ge \lambda),
\end{align}
where we have used Theorem \ref{Zhai's thm} in (\ref{eq1-lower-minimum})
and (\ref{eq3-lower-minimum}),
and the symmetry of the GFF in (\ref{eq2-lower-minimum}).
Take $\lambda = 
N^{2 - 2(\sqrt{\theta} - \sqrt{\theta - 2 \eta \sqrt{\theta}})^2 - \varepsilon}.$
By the lower bound of Theorem \ref{maximum}(i), the probability in 
(\ref{eq3-lower-minimum})
goes to $1$ as $N \to \infty$.
~~$\Box$ \\

Next, we prove the upper bound of Theorem \ref{maximum}(ii). \\
{\it Proof of the upper bound of Theorem \ref{maximum}(ii).}
Fix $\theta > 1$, $\varepsilon >0$, and $\eta >0.$
Set
$$\alpha_1 := 1 - 1/\sqrt{\log N},~~ 
\alpha_2 := 1 - 2/\sqrt{\log N} + K/(\log N)^{3/4},
~~m_N := \lceil 2 \theta (\log N)^{3/2} \rceil,$$
$$\lambda_N := (1 + K/2(\log N)^{1/4}) \frac{2}{\pi} N^2 
\log (N^{\alpha_1}/N^{\alpha_2}) m_N,$$
where $K$ is a sufficiently large positive constant.
We have for all $\lambda > 0$,
\begin{align}
&~~~P_0 ( |\{x \in \mathbb{Z}_N^2 : 
\frac{L_{\tau_{t_{\theta}}}^N (x) - t_{\theta}}
{\sqrt{2t_{\theta}}} \le -\eta \cdot 2\sqrt{2/\pi} \log N \} |
\ge \lambda) \notag \\
&\le P_0 ( |\{x \in \mathbb{Z}_N^2 : 
\sum_{j = 1}^{m_N} L_x^{(j)} [N^{\alpha_2}, N^{\alpha_1}] \le 
(\theta - 2 \eta \sqrt{\theta}) 4/\pi (\log N)^2 \}| \ge \lambda) 
\label{eq1-upper-local-time}\\
&+ P_0 (\sum_{j=0}^{m_N} \tau_x^{(j)} [N^{\alpha_2}, N^{\alpha_1}] \ge \lambda_N
~\text{for some}~x \in \mathbb{Z}_N^2 ) 
\label{eq2-upper-local-time}\\
&+ P_0 (\tau_{t_{\theta}} < \lambda_N). \label{eq3-upper-local-time}
\end{align}
By Lemma \ref{excursion-length} and 
\ref{concentration-estimate-inverse-local-time},
(\ref{eq2-upper-local-time}) and (\ref{eq3-upper-local-time}) 
go to $0$ as $N \to \infty$.
In analogy to the proof of Proposition \ref{localtime_excursion}, by Lemma \ref{green-function} and 
\ref{moment-estimate-local-time}, we have for all 
$x \in \mathbb{Z}_N^2$,
$$P_0 ( 
\sum_{j = 1}^{m_N} L_x^{(j)} [N^{\alpha_2}, N^{\alpha_1}] \le 
(\theta - 2 \eta \sqrt{\theta}) 4/\pi (\log N)^2)
\le N^{-2(\sqrt{\theta} - \sqrt{\theta - 2 \eta \sqrt{\theta}})^2 + o(1_N)}.$$
Therefore,
by taking $\lambda = N^{2-2(\sqrt{\theta} - \sqrt{\theta - 2 \eta \sqrt{\theta}})^2
+ \varepsilon}$, we can show that (\ref{eq1-upper-local-time}) goes to $0$ as
$N \to \infty.$ $\Box$ \\\\
{\bf Acknowledgements.} \\
The author would like to thank Professor Biskup for suggesting the problem
and stimulating discussions,
and thank Professor Zeitouni for 
fruitful comments which leads to
a significant improvement of an earlier work.
He would like to thank Professor Kumagai for variable comments
and encouragement.

\end{document}